\documentclass[12pt,psamsfonts]{amsart}
\usepackage{amsmath}
\usepackage{amsthm}
\usepackage{amscd}
\usepackage{amsfonts}
\usepackage{amsbsy}
\usepackage{latexsym}
\usepackage[psamsfonts]{amssymb}
\usepackage{graphicx}

\newtheorem{theorem}{Theorem}

\newtheorem{lemma}[theorem]{Lemma}

\newtheorem{remark}[theorem]{Remark}

\begin{document}

\title{Planar polynomial vector fields having a polynomial
first integral can be obtained from linear systems}
\author[B. Garc\'{\i}a, H. Giacomini and J.S. P\'{e}rez del R\'{\i}o]
{Belen Garc\'{\i}a$^{1}$, H\'{e}ctor Giacomini$^{2}$ and Jes\'{u}s
S. P\'{e}rez del R\'{\i}o$^{1}$}

\address{${^1}$ Departamento de Matem\'aticas, Universidad de Oviedo. Avda Calvo Sotelo, s/n., 33007, Oviedo, Spain.}

\email{belen.garcia@uniovi.es, jspr@uniovi.es}

\address{${^2}$ Laboratoire de Math\'{e}matiques et Physique Th\'{e}orique.
                                                    Facult\'{e} de Sciences et Techniques. Universit\'{e} Francois Rabelais de Tours.
                                                    Parc de Grandmont, 37200 Tours-France.
}

\email{Hector.Giacomini@lmpt.univ-tours.fr}

\thanks{The first and third authors are partially supported by a MEC/FEDER grant number
MTM2008-06065.} \subjclass{Primary 34C05, 34A34, 34C14}

\keywords{polynomial vector field, polynomial first integral,
linearization}

\begin{abstract}
We consider in this work planar polynomial differential systems having a polynomial first integral. We prove that these systems can be obtained from a linear system through a polynomial change of variables.

\end{abstract}

\maketitle \pagestyle{myheadings} \markboth{B. Garc\'{\i}a, H.
Giacomini and J. S. P\'{e}rez del R\'{\i}o}{Linearization of
Polynomial Vector Fields\ldots}

\section{Introduction and main results}

We consider planar vector fields of the form
\begin{equation}
x^{\prime }=P(x,y),\text{ \ \ \ \ \ }y^{\prime }=Q(x,y),  \label{uno}
\end{equation}%
where $P$ and $Q$ $\in \mathbb{R}[x,y]$ are coprime polynomials, that is, there is no non-constant polynomial which divides $P$ and $Q$. Here the prime indicates a derivative with respect the real independent variable $t$ and $\mathbb{R}[x,y]$ is the ring of all polynomials in the variables $x$ and $y$ with coefficients in $\mathbb{R}$. Let $m$ be the maximun degree of the polynomials  $P$ and $Q$, we say that system (\ref{uno}) is of degree $m$. We denote  $X=(P,Q)$  the vector field associated
to (\ref{uno}).

We say that system (\ref{uno}) has a polynomial first integral $H$
$\in \mathbb{R}[x,y]$ if and only if $\dfrac{\partial H}{\partial
x}P+\dfrac{\partial H}{\partial y}
Q\equiv 0$ on $\mathbb{R}^2$.\\
The search of first integrals is a classical tool for classifying
all trajectories of a polynomial system. A polynomial first
integral $H$ of (\ref{uno}) is called \textit{minimal} if $\deg
(H)\leq \deg (\tilde{H})$, for any polynomial first integral
$\tilde{H}$ of ( \ref{uno}). A particular class of systems having
a polynomial first integral is that of polynomial
\textit{Hamiltonian} systems, that is, systems of the form
$P=\dfrac{\partial H}{\partial y}$ , $Q=-
\dfrac{\partial H}{\partial x}$, where $H\in \mathbb{R}[x,y]$.\\
We say that an algebraic curve $ f(x,y)=0$, with $f(x,y) \in \mathbb{C}[x,y]$, is an invariant curve of
(\ref{uno}) if and
only if there exists a polynomial $K\in \mathbb{C}[x,y]$ such that $\dfrac{\partial f}{\partial x}P+\dfrac{\partial f}{\partial y}Q\equiv Kf$. Here $\mathbb{C}[x,y]$ is the ring of all polynomials in the variables $x$ and $y$ with coefficients in $\mathbb{C}$. The polynomial $K$ is called the cofactor associated to the invariant algebraic curve $f$.

If the polynomials $P$ and $Q$ are coprime then the degree of
a minimal polynomial first integral is greater than the degree of the system
(see \cite{LZ}).

We recall that $R(x,y)$  is an \textit{integrating factor} of (\ref{uno}) if it is a function of class $C^1$  in some open set $U$ of $\mathbb{R}^2$ and satisfies the partial differential equation
\begin{equation*}
\frac{\partial (RP)}{\partial x}=-\frac{\partial (RQ)}{\partial y}.
\end{equation*}%
In this case the system has a first integral given by
\begin{equation*}
H(x,y)=\int R(x,y)P(x,y)dy+f(x),
\end{equation*}%
where $f(x)$ must be choosen in such a way that the equality
$\dfrac{\partial H}{\partial x}
(x,y)=-R(x,y)Q(x,y)$ holds. Obviously, if $R$ is a polynomial then also $%
H$ is a polynomial and it is know (see \cite{LZ}\ ) that if system (\ref{uno}) has a minimal polynomial first
integral having degree greater than $m+1$ then there exists a polynomial
integrating factor with degree equal to $\deg (H)-m-1$.\\
We recall now the concept of remarkable value introduced by
Poincar\'{e} (see \cite{P}) in order
to obtain the classification of the polynomial vector fields of a given
degree $m$ having a rational first integral, and, in particular, a
polynomial first integral. In fact, he supposed that the degree of
the minimal polynomials that appear in the rational first integral have an upper bound that
depends of $m$.  Nevertheless this supposition was shown to be incorrect and this
 classification is still an open problem for rational first integrals, even for $m=2$. Some results in this direction were obtained in \cite{CL} and \cite{LO}. Very recently this problem was solved for the case of polynomial first integrals and $m=2$ (\cite{CGLPR}). They provide the necessary and sufficient conditions in order that a quadratic system have a polynomial first integral (see also \cite{ALV}). They also give the explicit polynomial first integrals in terms of the coefficients of a set of normal forms. A remarkable fact of such study is that quadratic systems may have polynomial first integrals of any arbitrary degree.
If the system (\ref{uno}) has a minimal polynomial first
integral $H$\ we say that $c\in \mathbb{C}$\ is a \textit{remarkable value}
for $H$ if $H+c$\ is not irreducible, that is, if there exist values $%
k_{1},..,k_{p}\in \mathbb{N}$ such that $%
H+c=u_{1}^{k_{1}}...u_{p}^{k_{p}}$\, where the $u_{i}$ are
irreducible polynomials on $\mathbb{C}[x,y]$ called \textit{remarkable factors
associated to }$c.$ Furthermore if some $k_{i}$ is greater than $1$
then the remarkable value is called \textit{critical} and the corresponding
factor $u_{i}$\ is called a \textit{critical remarkable factor}. The algebraic curves defined by the remarkable factors  appear to be very important in the construction of the phase portrait of the vector field, as it is shown in \cite{FL}.

In \cite{CGGL} the authors prove that the number of remarkable values
of a minimal polynomial first integral of a polynomial vector field is
finite.

If the polynomial first integral has no critical remarkable values
Javier Chavarriga proved that the system is Hamiltonian (see
\cite{FL}).

In this work we will obtain new results about the
polynomial vector fields having a  polynomial first integral. First of all we present a result that caracterizes the existence of exactly one critical remarkable value through the degree of the remarkable factors
that appear in the factorisation of the polinomial first integral, associated to a critical remarkable value. We can suppose, without loss of generality, that this remarkable value is $c=0$.\\

\textbf{Theorem A}. \textit{A minimal polynomial first integral of system (%
\ref{uno}) of the form} $%
H(x,y)=u_{1}(x,y)^{k_{1}}....u_{p}(x,y)^{k_{p}}$,\textit{
where the $u_{i}(x,y)$ are irreducible polynomials in $\mathbb{C}[x,y]$ and the $k_{i} \in \mathbb{N}$ with some  $k_{i}>1$, has exactly one critical remarkable value
if and only if $\deg (u_{1})+...+\deg (u_{p})=m+1$}.\\\\
This theorem will  be proved in section two.

In the following theorem we obtain the general expression of the
polynomial vector field having a given polynomial first integral.
We give conditions for  the degree of this vector field to be the
smallest one and we prove that if these conditions are satisfied
the polynomial first integral has exactly one critical remarkable
value. We say that the set of irreducible algebraic curves
$u_{i}(x,y)=0$ in $\mathbb{C}[x,y]$, with $i=1,...,p$ satisfy
\textit{the Christopher-Zoladek conditions} (these conditions
expressed in the same or similar form appear, for instance in
\cite{C}, \cite{Z}) if:

(i) there are no points at which $u_{i}$ and its two first partial derivatives simultaneously
vanish,

(ii) the highest order term of $u_{i}$ has not repeated factors,

(iii) no more than two curves meet at each point in the finite plane and when two curves meet at a point they are not tangent,

(iv) no two curves have a common factor in their highest order terms.\\

\textbf{Theorem B}.

\textit{Consider a polynomial }$%
H(x,y)=u_{1}(x,y)^{k_{1}}...u_{p}(x,y)^{k_{p}}$\textit{,
with $p>1$, where the} $u_{i}(x,y)$ \textit{are irreducible polynomials in $\mathbb{C}[x,y]$ and the $k_{i} \in \mathbb{N}$. Then }$H$
\textit{is a polynomial first integral of the polynomial vector
field defined by }
\begin{equation}
X=\sum\limits_{l=1}^{p}k_{l}\left( \prod\limits_{\substack{ i=1
\\ i\neq l}} ^{p}u_{i}\right) H_{u_{l}},  \label{eso}
\end{equation}%
\textit{where } $H_{u_{l}}=(\frac{\partial u_{l}}{\partial y},
-\frac{\partial u_{l}}{\partial x})$ \textit{\ is the Hamiltonian
vector field associated to the polynomial} $u_{l}$. \textit{
Furthermore, if the algebraic curves }$u_{i}(x,y)=0$ \textit{\
satisfy the Christopher-Zoladek conditions, then (\ref{eso}) is
the polynomial system of smallest degree that has }$H$\textit{\ as
a polynomial first integral and  }$H$\textit{\ has exactly one
critical remarkable
value. }\\

The following remark shows the importance of the Christopher-Zoladek conditions.\\

\textbf{Remark C}. \textit{If the algebraic curves $u_{i}(x,y)=0$
do not satisfy the Christopher-Zoladek conditions, then it can
exist vector fields of degree $m <\sum\limits_{i=1}^{p}\deg
(u_{i})-1$ having  $H$
 as polynomial first integral. In this case, the polynomials $P$ and $Q$ of the vector field
defined by (\ref{eso}) are not coprime}.\\
This remark will be justified in section three.
Finally we enonce the more important result of this work: a non hamiltonian polynomial vector field
having a polynomial first integral can be obtained from
a linear system through a polynomial change of variables.\\

\textbf{Theorem D}.
\textit{If a non Hamiltonian polynomial system (\ref{uno}) has a
polynomial first integral then it can be obtained from a linear system
through a polynomial change of variables and a rational change of the independent variable .
}\\

In Section two we prove Theorem A and
in Section three we prove Theorems B and D.

\section{Invariant curves and remarkable critical values.}

Let $H$ be a minimal polynomial first integral of system
(\ref{uno})  having remarkable critical values $c_{1},...,c_{s}$
and let $u_{1},...,u_{r}$ be all the remarkable factors (critical
or not) associated to these remarkable critical
values. We define $R=\prod%
\limits_{i=1}^{r}u_{i}^{k_{i}-1}$ and
$V=\prod\limits_{i=1}^{r}u_{i}$, with the $k_{i} \in \mathbb{N}$.
The factors that contribute in the expression of $R$ are only the
critical remarkable factors with their corresponding exponents
$k_{i}>1$. In \cite{CGGL} the authors proved that if in
(\ref{uno}) the polynomials $P$ and $Q$ are coprime then $R$ is a
polynomial integrating factor. This fact implies that $V$ is a
polynomial inverse integrating factor because it is the quotient
between the first integral $\tilde{H}=\prod\limits_{i=1}^{s}(H+c_{i})$ and
the integrating factor $R$ (for the definition and properties of
the inverse integrating factor see for instance \cite{CGGL}).

On the other hand, in \cite{FL} it is proved that if $H_{R}$ is the
first integral associated to the polynomial integrating factor $R$, then $H_{R}$ is a minimal polynomial first integral. Then $H_{R}$  and $H$ can only differ in multiplicative and additive constants. It is also proved in \cite{FL} that
\begin{equation}
\deg (V)=(s-1)d+(m+1)s,  \label{for}
\end{equation}%
where $s$ is the number of critical remarkable values of $H_{R}$, $d=\deg (R)$
and $m$ is the degree of the polynomial vector field.

We will prove now Theorem A.\\

\textit{Proof of Theorem A}. If the system has a minimal polynomial first
integral $H$ with exactly one critical remarkable value $c=0$ and with $%
H(x,y)=u_{1}(x,y)^{k_{1}}....u_{p}(x,y)^{k_{p}}$, where the $u_{i}$ are the
irreducible factors, we have $V=\prod\limits_{i=1}^{p}u_{i}$. Then from (\ref{for}) we have that $\deg (V)=\deg
(u_{1})+...+\deg (u_{p})=m+1$. We will prove now the converse result. We suppose that the last equality holds. If there exist more critical remarkable values, then the
integrating factor must be of the form\\
\begin{equation*}
R(x,y)=u_{1}(x,y)^{k_{1}-1}...u_{p}(x,y)^{k_{p}-1}w(x,y),
\end{equation*}
where $w$ is the product of the remarkable factors associated to
the other remarkable critical values. Since $H$ has a critical
remarkable value, the system is non Hamiltonian, and then
deg$(H)>m+1$. Then we have, as it was proved in \cite{LZ}, that
deg$(H)=m+1+$deg$(R)$. Hence we can write
\begin{equation*}
\begin{array}{l}
\displaystyle
 k_{1}\deg(u_{1})+\ldots +k_{p}\deg (u_{p})=\\
\displaystyle
(k_{1}-1)\deg
(u_{1})+...+(k_{p}-1)\deg (u_{p})+\deg (w)+m+1,
\end{array}
\end{equation*}
and we deduce that
\begin{equation*}
\deg (u_{1})+...+\deg (u_{k})=m+1+deg(w)>m+1,
\end{equation*}
a contradiction. Then the proof is finished.
$ \blacksquare $

\section{Polynomial vector fields with polynomial first integral obtained from linear vector fields}
From the above results it is easy to prove the following Lemma.\\

\begin{lemma}
\label{otra} If system (\ref{uno}) has a polynomial first integral and it is not Hamiltonian, then there exists a minimal
polynomial first integral $H$ that can be expressed as $%
H(x,y)=u_{1}(x,y)^{k_{1}}....u_{p}(x,y)^{k_{p}}$, where $p>1$, the $u_{i}(x,y)$
are irreducible polynomials in $\mathbb{C}$ and at least one of the $k_{i}$ is greater than $1$.
\end{lemma}

\textit{Proof}. The existence of a minimal polynomial first integral having
at least one remarkable critical value is clear because the non existence of
critical remarkable values implies that the system is Hamiltonian (see \cite%
{FL}). If $c=0$ is one of these remarkables values and $p=1$ with $H=u_{1}^{k_{1}}$ and $%
k_{1}>1$, then $H$ is not minimal and therefore the Lemma is proved. $%
\blacksquare $

\begin{remark}
Since we consider real systems, if among the irreducible polynomial factors $u_{i}$ there
is one that is not real then its complex conjugate must also appear as a factor, with the same exponent. Then the product of $u_{i}$ and its complex conjugate is
a real polynomial factor of the first integral.
\end{remark}

Therefore, if system (\ref{uno}) has a polynomial first integral
and it is not Hamiltonian, there exists a minimal polynomial first
integral having at least two real irreducible factors.
In this case we can obtain a polynomial vector field that has this polynomial first integral in the form given in Theorem B.\\

\textit{Proof of Theorem B}. In order to prove the first part
of the theorem it is sufficient to verify that $\dfrac{\partial H}{\partial x}P+\dfrac{\partial H}{\partial y}%
Q\equiv 0$. To prove the second part of the theorem we take into account that the degree of the vector field (\ref%
{eso}) is, by construction, $m=\sum\limits_{i=1}^{p}\deg
(u_{i})-1.$ Suppose that there exists another vector field
$\tilde{X}$ of smaller degree that has the invariant algebraic
curves $u_{1}(x,y)=0,...,u_{p}(x,y)=0$. As these curves satisfy
the Cristopher-Zoladek conditions , from Theorem 1 of \cite{CLPZ}
(this theorem was previously stated in others papers without proof
(\cite{C}, \cite{KC})), the vector field $\tilde{X}$ must be
identically zero. Therefore the minimun degree of the polynomial
vector field that has the polynomial first integral $H$ is $m$.
Finally, from Theorem A we deduce that $H$ has exactly one
remarkable
value. $\blacksquare $\\

If the set of irreducible curves do not satisfy the
Christopher-Zoladek conditions, then it is possible to have a
vector field $X$ such that $\deg (X)<\sum\limits_{i=1}^{p}\deg
(u_{i})-1$, having $
H(x,y)=u_{1}(x,y)^{k_{1}}....u_{p}(x,y)^{k_{p}}$ as a polynomial
first
integral.\\
 In the next Lemma we show the relationship between two polynomial vector
fields having the same polynomial first integral.\\

\begin{lemma}
Let $X_{1}$ and $X_{2}$ be two polynomial vector fields with different degrees
$m_{1}$ and $ m_{2}$ where $m_{2}>m_{1}$. We suppose that at least one of these vector fields (say $X_{1}$) has components given by coprime polynomials. If $X_{1}$ and $X_{2}$ have the same polynomial first integral
$H$ then $X_{2}=GX_{1}$ where $G$ is a polynomial of degree
$m_{2}-m_{1}$.
\end{lemma}

\textit{Proof}. If we denote $X_{1}=(P_{1},Q_{1})$ and $X_{2}=(P_{2},Q_{2})$
we deduce that $\dfrac{\partial H}{\partial x}P_{1}+\dfrac{\partial H}{\partial y}%
Q_{1}\equiv 0$ and $\dfrac{\partial H}{\partial x}P_{2}+\dfrac{\partial H}{\partial y}%
Q_{2}\equiv 0$, that is, $
P_{1}/Q_{1}=P_{2}/Q_{2}$. Since $P_{1}$ and $Q_{1}$ have no common factors we
deduce that there exists a polynomial $G$ of degree $m_{2}-m_{1}$ such that $%
P_{2}=GP_{1}$ and $Q_{2}=GQ_{1}$, and the proof is completed. $\blacksquare $
\\\\
One inmediate consequence of Theorem B and the above Lemma is the
Remark C that is stated in the Introduction.\\
A important consequence of Theorem B is that each non Hamiltonian
polynomial system with a polynomial first integral can be
obtained through a polynomial change of variables from a linear
system having a saddle point at the origin.\\

\textit{Proof of Theorem D}. From Lemma \ref{otra} we can suppose
that the polynomial system $x^{\prime }=P(x,y),$ $y^{\prime
}=Q(x,y)$ has a minimal first integral given by $
H(x,y)=u_{1}(x,y)^{k_{1}}....u_{p}(x,y)^{k_{p}}$ with $p>1$ and at
least one of the $k_{i}>1$, i.e. we suppose that $c=0$ is a
critical remarkable value. Consider the linear system $\dot{u}=u$,
$\dot{v}=-v$, where the dot indicates a derivative with respect to
an independent variable $\tau$. We introduce now the change of
variables $(u,v)\rightarrow (x,y)$ defined by
$u=u_{1}^{k_{1}}(x,y)...u_{p-1}^{k_{p-1}}(x,y)$ and $v=u_{p}^{k_{p}}(x,y)$.\\
A direct calculation gives
\begin{equation*}
\begin{array}{cc}
K_{1}(x,y)\dot{x}+K_{2}(x,y)\dot{y} & =\prod
\limits_{l=1}^{p-1}u_{l}, \\
K_{3}(x,y)\dot{x}+K_{4}(x,y)\dot{y} & =-u_{p},
\end{array}
\end{equation*}
where
\begin{equation*}
\begin{array}{cc}
K_{1}(x,y)=\left[\sum\limits_{i=1}^{p-1}k_{i}(u_{i})_{x}\prod \limits_{j=1, j\neq
i}^{p-1}u_{j}\right] \;, \; K_{2}(x,y)=\left[\sum\limits_{i=1}^{p-1}k_{i}(u_{i})_{y}\prod\limits_{ j=1, j\neq
i}^{p-1}u_{j}\right] , \\\\
K_{3}(x,y)=k_{p}\left( u_{p}\right) _{x} \; \; , \; \; K_{4}(x,y)=k_{p}\left( u_{p}\right) _{y}.
\end{array}
\end{equation*}
Here $(u_{i})_{x}$ and $(u_{i})_{y}$ indicate the partial derivatives of $u_{i}$ with respect to $x$ and $y$, respectively.

 If we denote $D(x,y)=K_{1}(x,y)K_{4}(x,y)-K_{2}(x,y)K_{3}(x,y)$ we obtain in the
region where the polynomial $D(x,y)\neq 0$
\begin{equation*}
\begin{array}{cc}
\dot{x} & =\dfrac{K_{4}\prod\limits_{l=1}^{p-1}u_{l}+K_{2}u_{p}}{D(x,y)}, \\
\dot{y} & =\dfrac{-K_{1}u_{p}-K_{3}\prod\limits_{l=1}^{p-1}u_{l}}{D(x,y)}.
\end{array}
\end{equation*}
We have also
\begin{equation*}
K_{4}\prod\limits_{l=1}^{p-1}u_{l}+K_{2}u_{p}= \left[\sum\limits_{i=1}^{p}k_{i}(u_{i})_{y}\prod\limits_{ j=1, j\neq
i}^{p}u_{j}\right],
\end{equation*}
\begin{equation*}
-K_{1}u_{p}-K_{3}\prod\limits_{l=1}^{p-1}u_{l}=-\left[ \sum\limits_{i=1}^{p}k_{i}(u_{i})_{x}\prod\limits_{ j=1, j\neq i}^{p}u_{j}\right].
\end{equation*}

The vector field $(\left[\sum\limits_{i=1}^{p}k_{i}(u_{i})_{y}\prod\limits_{ j=1, j\neq
i}^{p}u_{j}\right],-\left[ \sum\limits_{i=1}^{p}k_{i}(u_{i})_{x}\prod\limits_{ j=1, j\neq i}^{p}u_{j}\right])$ is the vector field that appears in (2) and then it has the polynomial first integral H. Therefore, as  the polynomials $P$ and $Q$ are coprime, by Lemma 3 there exists a polynomial $G$ such that
\begin{equation*}
\left[\sum\limits_{i=1}^{p}k_{i}(u_{i})_{y}\prod\limits_{ j=1, j\neq
i}^{p}u_{j}\right]=GP \; , \; -\left[ \sum\limits_{i=1}^{p}k_{i}(u_{i})_{x}\prod\limits_{ j=1, j\neq i}^{p}u_{j}\right]=GQ.
\end{equation*}

 Then we have
\begin{equation*}
\begin{array}{cc}
\dot{x} & =\dfrac{G(x,y)P(x,y)}{D(x,y)}, \\
\dot{y} & =\dfrac{G(x,y)Q(x,y)}{D(x,y)}.
\end{array}
\end{equation*}
Hence through a rational change of the independent variable $\tau\rightarrow t$ defined by $d\tau=\frac{D(x,y)}{G(x,y)}dt$ in the region where $D(x,y)\neq 0$ and $G(x,y)\neq 0$, we obtain the original system $x^{\prime }=P(x,y)$, $y^{\prime }=Q(x,y)$. $\blacksquare $ \\\\
Obviously, the change of variables that transforms the linear vector field $(u,-v)$ into the vector field $(P,Q)$ is not a diffeomorphism in all $\mathbb{R}^2$. The qualitative behaviour of the orbits is completely changed by this change of variables. It simply enables us to obtain through algebraic calculations all non-hamiltonian polynomial vector fields with a polynomial first integral from a linear hamiltonian vector field with a saddle point at the origin.


\begin{thebibliography}{9}
\bibitem{ALV} J. Artes, J. Llibre and N. Vulpe, \textit{Quadratic systems with a polynomial first integral: A complete classification in the coefficients space $\mathbb{R}^{12}$ }, J. Differential Equations \textbf{246} (2009), 3535-3558.

\bibitem{CL} L. Cairo and J. Llibre, \textit{Phase portraits of quadratic polynomial vector fields having a rational first integral of degree 2}, Nonlinear Analysis \textbf{67} (2007), 327-348.

\bibitem{CGLPR} J. Chavarriga, B. Garc\'{\i}a, J. Llibre, J. S. P\'{e}rez
del R\'{\i}o and J. A. Rodr\'{\i}guez, \textit{Polynomial first integrals of
quadratic vector fields, }J. Differential Equations \textbf{230} (2006),
393--421.

\bibitem{CGGL} J. Chavarriga, H. Giacomini, J. Gin\'{e} and J. Llibre,
\textit{Darboux integrability and the inverse integrating factor, }J.
Differential Equations \textbf{194} (2003), 116--139.

\bibitem{C} C. Christopher, \textit{ Invariant algebraic
curves and conditions for a centre , }Proc. Roy. Soc.
Edinburgh \textbf{124A }(1994), 1209-1229.

\bibitem{CLPZ} C. Christopher, J. Llibre, C. Pantazi and X. Zhang,
\textit{Darboux integrability and invariant algebraic curves for
planar polynomial systems}, J. Phys. A: Math. Gen. \textbf{35}
(2002), 2457-2476.

\bibitem{FL} A. Ferragut and J. Llibre, \textit{On the 
remarkable values of the rational first integrals of polynomial
vector fields, }J. Differential Equations \textbf{241} (2007),
399-417.

\bibitem{KC} R. E. Kooij and C. J. Christopher, \textit{
Algebraic invariant curves and the integrability of polynomial
systems , }Appl. Math. Lett. \textbf{6, }(1993),
51-53.

\bibitem{LO} J. Llibre and R. Oliveira, \textit{Phase portraits of quadratic polynomial vector fields having a rational first integral of degree 3},
Nonlinear Analysis \textbf{70} (2009), 3549-3560.

\bibitem{LZ} J. Llibre and X. Zhang, \textit{Polynomial first integrals of
quadratic systems}, Rocky Mountain J. of Math. \textbf{31} (2001), 1317-1371.

\bibitem{P} H. Poincar\'{e}, \textit{Sur l'int\'{e}gration des \'{e}quations
diff\'{e}rentielles du premier ordre et du premier degr\'{e} I and II},
Rendiconti del Circolo Matematico di Palermo \textbf{5} (1891), 161--191;
\textbf{11} (1897), 193--239.

\bibitem{Z} H. Zoladek, \textit{On algebraic solutions of
algebraic Pfaff equations , }Studia Mathematica \textbf{%
114 }(1995), 117-126.
\end{thebibliography}
\end{document}